\documentclass{agtart_a}
\pdfoutput=1
\usepackage{graphicx}

\title{Lattices acting on right-angled buildings}

\author{Anne Thomas}
\givenname{Anne}
\surname{Thomas}
\address{Department of Mathematics\\
University of Chicago\\\newline
5734 South University Ave\\
Chicago IL 60637\\
USA}
\email{athomas@math.uchicago.edu}
\urladdr{}

\volumenumber{6}
\issuenumber{}
\publicationyear{2006}
\papernumber{42}
\startpage{1215}
\endpage{1238}

\doi{}
\MR{}
\Zbl{}


\received{3 March 2006}
\revised{29 June 2006}
\accepted{30 June 2006}
\published{}
\publishedonline{7 September 2006}
\proposed{}
\seconded{}
\corresponding{}
\editor{}
\version{}

\arxivreference{math.GR/0508385}

\let\xysavmatrix\xymatrix
\def\xymatrix{\disablesubscriptcorrection\xysavmatrix}
\AtBeginDocument{\let\bar\wbar\let\tilde\wtilde\let\hat\what}

\newcommand{\N}{\ensuremath{\mathbb{N}}}

\newcommand{\Hyp}{\ensuremath{\mathbb{H}}}
\newcommand{\cA}{\ensuremath{\mathcal{A}}}
\newcommand{\cB}{\ensuremath{\mathcal{B}}}
\newcommand{\cG}{\ensuremath{\mathcal{G}}}

\newcommand{\cC}{\ensuremath{\mathcal{C}}}
\newcommand{\cS}{\ensuremath{\mathcal{S}}}
\newcommand{\cV}{\ensuremath{\mathcal{V}}}

\newcommand{\Aut}{\ensuremath{\operatorname{Aut}}}
\newcommand{\Fix}{\ensuremath{\operatorname{Fix}}}
\newcommand{\Vol}{\ensuremath{\operatorname{Vol}}}

\newcommand{\Ad}{\ensuremath{\operatorname{Ad}}}

\newcommand{\St}{\ensuremath{\operatorname{St}}}

\newcommand{\ex}{\ensuremath{\operatorname{ex}}}

\newcommand{\quot}{\backslash \! \backslash}

\newcommand{\G}{\Gamma}

\newcommand{\bs}{\backslash}

\makeatletter
\def\cnewtheorem#1[#2]#3{\newtheorem{#1}{#3}[section]
\expandafter\let\csname c@#1\endcsname\c@theorem}
\makeatother

\newtheorem{theorem}{Theorem}[section]
\cnewtheorem{lemma}[theorem]{Lemma}
\cnewtheorem{corollary}[theorem]{Corollary}
\cnewtheorem{proposition}[theorem]{Proposition}
\newtheorem{maintheorem}{Main Theorem}
\makeautorefname{maintheorem}{Main Theorem}

\newtheorem{functortheorem}{Functor Theorem}
\makeautorefname{functortheorem}{Functor Theorem}

\numberwithin{equation}{section}

\begin{document}

\begin{asciiabstract}
Let X be a right-angled building.  We show that the lattices in
Aut(X) share many properties with tree lattices.  For example, we
characterise the set of covolumes of uniform and of nonuniform
lattices in Aut(X), and show that the group Aut(X) admits an
infinite ascending tower of uniform and of nonuniform lattices.
These results are proved by constructing a functor from graphs of
groups to complexes of groups.
\end{asciiabstract}

\begin{abstract}
Let $X$ be a right-angled building.  We show that the lattices in
$\mathrm{Aut}(X)$ share many properties with tree lattices.  For example,
we characterise the set of covolumes of uniform and of nonuniform
lattices in $\mathrm{Aut}(X)$, and show that the group $\mathrm{Aut}(X)$
admits an infinite ascending tower of uniform and of nonuniform lattices.
These results are proved by constructing a functor from graphs of groups
to complexes of groups.
\end{abstract}

\maketitle

%%%%%%%%%%%%%%%%%%%%%%%%%%%%%%%%%%%%%%%%%%%%%%%%%%%%%%%%%%%%%%%%%%%%
\section*{Introduction}\label{s:intro}
%%%%%%%%%%%%%%%%%%%%%%%%%%%%%%%%%%%%%%%%%%%%%%%%%%%%%%%%%%%%%%%%%%%%

Let $G$ be a locally compact topological group, with suitably
normalised Haar measure $\mu$.  A discrete subgroup $\G \leq G$ is a
\emph{lattice} if the covolume $\mu(\G\bs G)$ is finite, and a
\emph{uniform lattice} if $\G \bs G$ is compact.  A \emph{tower} of
lattices is a strictly increasing infinite chain of subgroups\[ \G_1
< \G_2 < \cdots < \G_i < \cdots
\] such that each $\G_i$ is a lattice in $G$.  Two basic questions are:

\begin{enumerate}
\item  What are the possible covolumes of lattices in
$G$? \item  Does $G$ admit a tower of lattices?
\end{enumerate}

These questions have been studied for many $G$. For example, if $G$
is a non-compact simple real Lie group, then the covolumes of
lattices in $G$ are bounded away from 0, and in most cases the set
of lattice covolumes is discrete (see Lubotzky~\cite{l1:tl} and the
references therein). A strong finiteness result for lattices in
semisimple groups is that of Borel~\cite{b1:ds}, which implies that
for $G$ a higher-rank algebraic group over a local field, and any $c
> 0$, there are only finitely many lattices in $G$ with covolume less than $c$.

Nonclassical cases arise from the fact that if $G$ is the
automorphism group of a locally finite polyhedral complex, then $G$
is naturally a locally compact group (see
\fullref{ss:lattcovols}). Lattices in the automorphism groups of
trees are treated by Bass and Lubotzky~\cite{bl1:tl}. Lattices acting on
a product of trees have been studied by, for example,
Burger--Mozes~\cite{bm1:lpt}.

We consider lattices in the automorphism groups of certain
buildings.  Let $P$ be a compact, convex polyhedron in $\Hyp^n$ with
all dihedral angles $\frac{\pi}{2}$, and let $(W,I)$ be the
right-angled Coxeter system generated by reflections in the
$(n-1)$--dimensional faces of $P$.  A \emph{right-angled hyperbolic
building} of type $(W,I)$ is a polyhedral complex $X$, equipped with
a maximal family of subcomplexes called \emph{apartments}, each
isometric to the tesselation of $\Hyp^n$ by copies of $P$ (called
\emph{chambers}), satisfying the usual axioms for a building.  A
right-angled hyperbolic building has nonpositive curvature. It is
\emph{locally} the product of trees, but is \emph{not} globally a
product. In dimension 2, where right-angled hyperbolic buildings are
sometimes known as ``Bourdon buildings", $P$ is a regular
right-angled hyperbolic polygon, and the link of each vertex of $X$
is a complete bipartite graph.

The dimension of a right-angled building with apartments isometric
to $\Hyp^n$ is at most 4, and this bound is
sharp (see Potyagailo--Vinberg~\cite{pv1:rrg}, and
Januszkiewicz--\'Swi\c{a}tkowski~\cite{js1:hcg}). However, given any
right-angled Coxeter system $(W,I)$, there exist right-angled
buildings with apartments isometric to the Davis complex for
$(W,I)$; thus right-angled buildings may be constructed in
arbitrarily high dimensions~\cite{js1:hcg}. In
\fullref{ss:hypbuildings} we give precise definitions and recall
the classification of right-angled buildings.  For work on
hyperbolic buildings see, for example, Gaboriau--Paulin~\cite{gp1:ih}
and its references.

Our results for lattices acting on right-angled buildings are as
follows.

\begin{maintheorem}
\label{thm:main}
Let $(W,I)$ be a right-angled Coxeter system and let $\{ q_i
\}_{i \in I}$ be a family of positive integers, $q_i \geq 2$, with
$q_{i_1}=\max_i\{q_i\}$.  Suppose $X$ is the (unique) building of type
$(W,I)$, such that each $\{i\}$--residue of $X$ has cardinality $q_i$.
Let $G=\Aut(X)$.
{\begin{enumerate} \item \emph{\textbf{Uniform covolumes}}\label{item:uniform}
\begin{enumerate}\item \label{i:unifcovol} If $m_{i_1,j}=\infty$ for
some $j \in I$, then the set of covolumes of uniform lattices in $G$
is
\[ \mathcal{V}_u(G) = \bigl\{ \tfrac{a}{b} ~\big|~
\mbox{prime divisors of $b$ are strictly less than $q_{i_1}$}
\bigr\}\]
\item \label{i:countunif}  If $m_{i_1,j}=\infty$ for
some $j \in I$, and $q_{i_1}> 2$, then for each $v \in
\mathcal{V}_u(G)$ there is a countably infinite number of
nonconjugate uniform lattices of covolume $v$.\end{enumerate}
\item \emph{\textbf{Nonuniform covolumes}}\label{item:nonuniform}
\begin{enumerate}
\item \label{i:nonunifcovol} If there exist $i, j \in I$ such that $m_{i,j}=\infty $ and $q_i > 2$, then the set of covolumes of nonuniform lattices is $(0,\infty)$.
 \item \label{i:nonunifcount} If there exist $i, j \in I$ such that
 $m_{i,j}=\infty$ and $q_i, q_j > 2$, then for every $v
> 0$ there exist uncountably many commensurability classes of
nonuniform lattices of covolume $v$.
\end{enumerate}
\item \emph{\textbf{Towers}}\label{item:towers}
\begin{enumerate}
\item \label{i:uniftower} If there exist $i, j \in I$ such that
$m_{i,j}=\infty $ and $q_i > 2$, then there exists a tower of uniform
lattices in $G$.
\item \label{i:nonuniftower} If there exist $i, j \in I$ such that $m_{i,j}=\infty $ and $q_i > 2$, then there exists a tower
of nonuniform lattices in $G$.
\item \label{i:gs} If there exist $i,
j \in I$ such that $m_{i,j}=\infty $ and $q_i$ is composite, then
there exists a tower of uniform lattices in $G$ such that the
quotient by each lattice is isometric to a chamber of
$X$.\end{enumerate}
\end{enumerate}}\end{maintheorem}

\noindent We observe that if $(W,I)$ is generated by reflections in
the faces of a right-angled hyperbolic polyhedron, then for each $i
\in I$, there is a $j \in I$ so that $m_{i,j} = \infty$.

Although the right-angled building $X$ is a higher-dimensional
object, the properties of lattices in $\Aut(X)$ listed in the
\fullref{thm:main} are quite different to those of lattices in higher-rank
algebraic groups. For example, the results cited above imply that
semisimple groups do not admit lattices of arbitrarily small
covolume, thus do not admit towers of (any) lattices.  The theorem
of Borel also implies that any $v
> 0$ is the covolume of at most finitely many lattices.  We note
too that there exist nonpositively curved polyhedral complexes whose
automorphism groups, while nondiscrete, do \emph{not} admit towers.
For example the Bruhat--Tits building for $SL_3(\Q_p)$ has
automorphism group containing $SL_3(\Q_p)$ as a finite index
subgroup.

In contrast, by comparing the \fullref{thm:main} with results for trees,
we see that the lattices in $G=\Aut(X)$ share many properties with
tree lattices.

\textbf{Uniform covolumes}\qua  Rosenberg
\cite[Proposition~9.1.2]{ros1:tctl} showed that the set of
covolumes of uniform lattices acting on the $m$--regular tree is
\[\bigl\{\tfrac{a}{b} ~\big|~ \mbox{prime divisors of $b$ are
  strictly less than $m$}\bigr\}\]
and found a similar result for
biregular trees \cite[Theorem~9.2.1]{ros1:tctl}.

Note that if all $q_i = 2$ then the set of uniform covolumes is
discrete. However, if some $q_i> 2$ then the set of (uniform)
covolumes is dense in $(0,\infty)$, and, in particular, $G$ admits
(uniform) lattices of arbitrarily small covolume, as is the case for
trees.

Our proof of the counting result~\eqref{i:countunif} will show that
this property holds for regular and biregular trees as well (see
\fullref{p:treecount} below). This generalises 
a result of Bass and Kulkarni \cite[Theorem~7.1(b)]{bk1:utl} and answers
a question put to us by A Lubotzky.

\textbf{Nonuniform covolumes}\qua Bass--Lubotzky
\cite[Theorem~4.3]{bl1:tl} showed that for the $m$--regular tree, $m \geq
3$, every $v > 0$ is the covolume of some nonuniform lattice. Rosenberg
\cite[Theorem~8.2.2]{ros1:tctl} extended this result to
biregular trees.

Farb--Hruska~\cite{fh1:ci} constructed commensurability invariants
for lattices acting on the $(m,n)$--biregular tree, for $m,n \geq
3$, and used these to show that for every $v > 0$ there are
uncountably many commensurability classes of nonuniform tree
lattices \cite[Corollary~1.2]{fh1:ci}.  Farb--Hruska have
shown~\eqref{i:nonunifcovol}
and~\eqref{i:nonunifcount} for some
$2$--dimensional right-angled hyperbolic buildings.

\textbf{Towers}\qua Rosenberg showed \cite[Theorem~3.3.1]{ros1:tctl}
that if $T$ is a tree such that $\Aut(T)$ is nondiscrete and admits a
uniform lattice, then $\Aut(T)$ admits a tower of uniform lattices, and
Carbone--Rosenberg showed \cite[Theorem~5.4]{cr1:it} that, with one
exception, if $\Aut(T)$ admits a nonuniform lattice then it admits a
tower of nonuniform lattices.

Part~\eqref{i:gs} of the \fullref{thm:main} addresses the following finer
question about towers. Let $X$ be a polyhedral complex and $G =
\Aut(X)$. A subgroup of $G$ is \emph{homogeneous} if it acts
transitively on the cells of maximum dimension in $X$.  Does $G$
admit a tower of homogeneous lattices?

When $X$ is the $3$--regular tree, a deep theorem of
Goldschmidt~\cite{g1:atg} implies that $G$ does not admit such a
tower, since $G$ has only finitely many conjugacy classes of
homogeneous lattices. The Goldschmidt--Sims conjecture
(see Glasner~\cite{g1:gs}), which remains open, is that if $X$ is the
$(p,q)$--biregular tree, where $p$ and $q$ are prime, then there are
only finitely many conjugacy classes of homogeneous lattices in $G$.
If $X$ is the product of two trees of prime valence,
Glasner~\cite{g1:gs} has shown that there are only finitely many
conjugacy classes of (irreducible) homogeneous lattices.  We do not
know if $G$ of the \fullref{thm:main} admits a homogeneous tower when all
$q_i$ are prime.  We note that, by the \fullref{thm:functor} below, if
such a $G$ does have only finitely many conjugacy classes of
homogeneous lattices, then the Goldschmidt--Sims conjecture holds.

The key fact used to prove the \fullref{thm:main} is
that if a tree is ``nicely embedded" in a right-angled building,
then any group of automorphisms of the tree may be extended to a
group of automorphisms of the building.  We express this by means of
a functor from graphs of groups to complexes of groups:

\begin{functortheorem}
Let $(W,I)$ and $X$ be as in the \fullref{thm:main}.  For each $i,
j \in I$ such that $m_{i,j}=\infty$, let $T$ be the $(q_{i},
q_{j})$--biregular tree.  Then there is a functor from the category
of graphs of groups with universal covering tree $T$ to the category
of complexes of groups with universal cover $X$, provided either
that $X$ is ``sufficiently symmetric", or that we restrict to
$2$--colourable graphs.  Moreover, this functor takes coverings to
coverings, and faithful graphs of groups to faithful complexes of
groups.
\end{functortheorem}

We give a precise statement and proof of the \fullref{thm:functor} in
\fullref{ss:functorthm}, and then apply it to prove the \fullref{thm:main} in
\fullref{s:proofmain}. An example of a sufficiently symmetric building
$X$ is one with regular chambers, for instance if $\dim(X)=2$, and all
$q_{i}$ equal.  We note that there are many variants of the
\fullref{thm:functor}, allowing promotion of known constructions of tree lattices to
more general complexes of groups. For the theory of graphs of groups and
their morphisms, see \fullref{ss:graphsgroups} and Serre~\cite{s1:t},
Bass~\cite{b1:ctgg} and Bass--Lubotzky~\cite{bl1:tl}. The theory of
complexes of groups (see Bridson and Haefliger~\cite{bh1:ms}) is outlined
in \fullref{ss:cxsgroups}.

I thank Benson Farb for introducing me to this area, his enthusiasm
and guidance, and his many comments on earlier versions of this
paper.  I would also like to thank G. Christopher Hruska and
Fr\'ed\'eric Haglund for helpful discussions, and the anonymous
referee for many worthwhile suggestions.

%%%%%%%%%%%%%%%%%%%%%%%%%%%%%%%%%%%%%%%%%%%%%%%%%%%%%%%%%%%%%%%%%%%
%%%%%%%%%%%%%%%%%%%%%%%%%%%%%%%%%%%%%%%%%%%%%%%%%%%%%%%%%%%%%%%%%%%
\section{Background}\label{s:background}
%%%%%%%%%%%%%%%%%%%%%%%%%%%%%%%%%%%%%%%%%%%%%%%%%%%%%%%%%%%%%%%%%%%
%%%%%%%%%%%%%%%%%%%%%%%%%%%%%%%%%%%%%%%%%%%%%%%%%%%%%%%%%%%%%%%%%%%

\fullref{ss:lattcovols} gives the basic definitions for lattices and
describes a suitable normalisation of Haar measure for automorphism
groups of polyhedral complexes. In \fullref{ss:hypbuildings} we discuss
right-angled buildings. The key definitions for graphs of groups
are given in \fullref{ss:graphsgroups}, and those for complexes of
groups in \fullref{ss:cxsgroups}.  We refer the reader to Bridson and
Haefliger~\cite{bh1:ms} for generalities on polyhedral complexes.

%%%%%%%%%%%%%%%%%%%%%%%%%%%%%%%%%%%%%%%%%%%%%%%%%%%%%%%%%%%%%%
\subsection{Lattices and covolumes}\label{ss:lattcovols}
%%%%%%%%%%%%%%%%%%%%%%%%%%%%%%%%%%%%%%%%%%%%%%%%%%%%%%%%%%%%%%

Let $G$ be a locally compact topological group with left-invariant
Haar measure $\mu$.  A discrete subgroup $\Gamma \leq G$ is a
\emph{lattice} if the covolume $\mu(\Gamma\backslash G)$ is finite.
A lattice $\G$ is \emph{uniform} if  $\G \bs G$ is compact.  Let $S$
be a left $G$--set such that for every $s \in S$, the stabiliser
$G_s$ is compact and open. Then if $\G \leq G$ is discrete, the
stabilisers $\Gamma_s$ are finite. We define the
\emph{$S$--covolume} of $\Gamma$ by
\[ \Vol(\Gamma \quot S) :=  \sum_{s
\in \Gamma \bs S } \frac{1}{|\Gamma_s|} \,\leq\infty
\] The following theorem shows that Haar measure may
be normalised so that $\mu(\G\bs G)$ equals the $S$--covolume.

\begin{theorem}[Serre \cite{s:cgd}]\label{t:Scovolumes} Let $G$ be a
locally compact topological group acting on a set $S$ with compact
open stabilisers and a finite quotient $G\backslash S$.  Suppose
further that $G$ admits at least one lattice.  Then there is a
normalisation of the Haar measure $\mu$, depending only on the
choice of $G$--set $S$, such that for each discrete subgroup
$\Gamma$ of $G$ we have $\mu(\Gamma \bs G) = \Vol(\Gamma \quot S)$.
\end{theorem}

Let $X$ be a connected, locally finite $n$--dimensional polyhedral
complex, with $X_n$ the set of $n$--dimensional cells of $X$.
Let $\Aut(X)$ be the group of cellular isometries of $X$.  A subgroup
of $\Aut(X)$ is said to act \emph{without inversions} if its elements
fix pointwise each cell that they preserve. The group $G = \Aut(X)$
is a locally compact topological group, with a neighbourhood basis
of the identity consisting of automorphisms fixing larger and larger
combinatorial balls.  By the same arguments as for tree lattices~(see
Bass and Lubotzky~\cite[Chapter~1]{bl1:tl}), it can be shown that if
$G\bs X$ is finite, then a discrete subgroup $\G \leq G$ is a lattice
if and only if its $X_n$--covolume converges, and $\G$ is uniform if and
only if this sum has finitely many terms. Using \fullref{t:Scovolumes},
we now normalise the Haar measure $\mu$ on $G=\Aut(X)$ so that for
all lattices $\G \leq G$, the covolume of $\G$ is
\[ \mu(\G \bs G) = \Vol(\G \quot X_n)\]

%%%%%%%%%%%%%%%%%%%%%%%%%%%%%%%%%%%%%%%%%%%%%%%%%%%%%%%%%%%%%%%%%%%
\subsection{Right-angled buildings}\label{ss:hypbuildings}
%%%%%%%%%%%%%%%%%%%%%%%%%%%%%%%%%%%%%%%%%%%%%%%%%%%%%%%%%%%%%%%%%%%

Let $(W,I)$ be a right-angled Coxeter system.  Let $N$ be the finite
nerve of $(W,I)$ and let $P'$ be the simplicial cone on $N'$ with
vertex $x_0$. We write $\mathcal{S}^f$ for the set of $J \subseteq
I$ such that the subgroup $W_J$ of $W$ generated by $J$ is finite.
By convention, $W_\phi = 1$, so the empty set $\phi$ is in $\cS^f$.
There is then a one-to-one correspondence between the vertices of
$P'$ and the types $J \in \cS^f$.  For each $i \in I$, the vertex of
$P'$ of type $\{i\}$ will be called an \emph{$i$--vertex}, and the
union of the simplices of $P'$ which contain the $i$--vertex but not
the cone point $x_0$ will be called an \emph{$i$--face}.

A \emph{right-angled building} of type $(W,I)$ is a polyhedral
complex $X$ equipped with a maximal family of subcomplexes, called
\emph{apartments}.  Each apartment is polyhedrally isometric to the
Davis complex for $(W,I)$, and the copies of $P'$ in $X$ are called
\emph{chambers}.  The apartments and chambers of $X$ satisfy the
usual axioms for a Bruhat--Tits building.  By abuse of notation, if
$\dim(X)=n$ then we write $X_n$ for the set of chambers of $X$.

Each vertex of a right-angled building $X$ has a type $J \in \cS^f$,
induced by the types of $P'$. For $i \in I$, an
\emph{$\{i\}$--residue} of $X$ is then the connected subcomplex
consisting of all chambers which meet a given $i$--face. Suppose
$q_i$ is the cardinality of each $\{i\}$--residue of $X$, that is,
the number of copies of $P'$ in each $\{i\}$--residue.  If $x$ is a
vertex of $X$, of type some $J \in \cS^f$ with $|J|=\dim(X)=n$, then
the link of $x$ in $X$ is the join of $n$ sets of points of
cardinalities respectively $q_j$, for $j \in J$.

If the Coxeter system $(W,I)$ is generated by reflections in the
faces of an $n$--dimensional right-angled hyperbolic polyhedron $P$,
then $P'$ is the barycentric subdivision of $P$.  The apartments of
$X$ are isometric to $\Hyp^n$, and we refer to $X$ as a hyperbolic
building.  For example, in Bourdon's building $I_{p,q}$, $P$ is a
regular right-angled hyperbolic $p$--gon (see Bourdon~\cite{b1:ih}). The
link of each vertex of $X$ is the complete bipartite graph
$K_{q,q}$, which may be thought of as the join of 2 sets of $q$
points, and each $\{i\}$--residue consists of $q$ copies of $P$,
glued together along a common edge.

The following result classifies right-angled buildings.

\begin{theorem}[{{Haglund--Paulin \cite[Proposition~1.2]{hp1:cai}}}]
\label{t:rabeu}
Let $(W,I)$ be a right-angled Coxeter
system and $\{q_i\}_{i \in I}$ a family of positive integers ($q_i
\geq 2$).  Then, up to isometry, there exists a unique building $X$
of type $(W,I)$, such that for each $i \in I$, the $\{i\}$--residue
of $X$ has cardinality $q_i$.
\end{theorem}

\noindent In the $2$--dimensional case, this result is due to
Bourdon~\cite{b1:ih}.  According to Hagland and Paulin~\cite{hp1:cai},
\fullref{t:rabeu} was proved by M Globus, and was known also to
M Davis, T Janusz\-kiewicz and J \'Swi\c{a}tkowski.

%%%%%%%%%%%%%%%%%%%%%%%%%%%%%%%%%%%%%%%%%%%%%%%%%%%%%%%%%%%%%%%%%%%%
\subsection{Graphs of groups}\label{ss:graphsgroups}
%%%%%%%%%%%%%%%%%%%%%%%%%%%%%%%%%%%%%%%%%%%%%%%%%%%%%%%%%%%%%%%%%%%%

We give only the definitions most relevant to the proof of the
\fullref{thm:functor} below.  See Serre~\cite{s1:t}, Bass--Lubotzky~\cite{bl1:tl}
and Bass~\cite{b1:ctgg} for more complete treatments.

A \emph{graph of groups} $\mathbb{A} = (A, \cA)$ consists of a
connected graph $A$, with vertices $V(A)$ and edges $E(A)$, together
with groups $\cA_v$ for each $v \in V(A)$ and $\cA_e =
\cA_{\overline{e}}$ for each $e \in E(A)$, and monomorphisms
$\alpha_e \co\cA_e \to \cA_{i(e)}$. See Bass~\cite{b1:ctgg} for the
definitions of the \emph{path group} $\pi(\mathbb{A})$, the
\emph{fundamental group of the graph of groups}
$\pi_1(\mathbb{A},v_0)$ and the \emph{universal covering tree}. A
graph of groups is \emph{faithful} if its fundamental group acts
faithfully on its universal covering tree.

Let $\mathbb{A} = (A, \cA)$ and $\mathbb{B} = (B, \cB)$ be graphs of
groups.  A \emph{morphism} $\phi\co\mathbb{A} \to \mathbb{B}$
consists of: \begin{enumerate}\item a morphism of graphs $f\co A \to
B$;\item homomorphisms of local groups $\phi_v\co \cA_v \to
\cB_{f(v)}$ and $\phi_e = \phi_{\overline{e}}\co\cA_e \to
\cB_{f(e)}$; and \item elements $\gamma_v \in \pi_1(\mathbb{B},
f(v))$ for each $v \in V(A)$, and $\gamma_e \in \pi(\mathbb{B})$ for
each $e \in E(A)$, such that if $v = i(e)$ then
\begin{enumerate}\item  $\delta_e := \gamma_v^{-1} \gamma_e \in
\cB_{f(v)}$; and\item $\phi_a\circ \alpha_e = \Ad(\delta_e)\circ
\alpha_{f(e)} \circ \phi_e$, where $\Ad(\delta_e)$ is conjugation by
$\delta_e$ in $\cB_{f(v)}$.
\end{enumerate}
\end{enumerate}
The morphism $\phi$ is a \emph{covering} if in
addition:\begin{enumerate} \item each $\phi_v$ and $\phi_e$ is
injective; and \item for each $v \in V(A)$ and $e' \in E(B)$ with
$i(e')=f(v)$, the map
\[  \coprod_{\substack{e \in f^{-1}(e') \\ i(e)=v}} \cA_v/\alpha_e(\cA_e)
\to \cB_{f(v)}/ \alpha_{e'}(\cB_{e'})\]
induced by $g \mapsto \phi_v(g) \delta_e$ is a bijection.
\end{enumerate}

%%%%%%%%%%%%%%%%%%%%%%%%%%%%%%%%%%%%%%%%%%%%%%%%%%%%%%%%%%%%%%%%%%%%
\subsection{Complexes of groups}\label{ss:cxsgroups}
%%%%%%%%%%%%%%%%%%%%%%%%%%%%%%%%%%%%%%%%%%%%%%%%%%%%%%%%%%%%%%%%%%%%

We sketch the theory of complexes of groups, due to Haefliger
\cite{bh1:ms}.  \fullref{ss:towers} outlines the construction to be used
for towers of lattices.

Throughout this section, if $X$ is a polyhedral complex, then $X'$
is the first barycentric subdivision of $X$. This is a simplicial
complex with vertices $V(X')$ and edges $E(X')$.  Each $a \in E(X')$
corresponds to cells $\tau \subset \sigma$ of $X$, and so may be
oriented from $\sigma$ to $\tau$.  We write $i(a)=\sigma$ and $t(a)
= \tau$. Two edges $a$ and $b$ of $X'$ are \emph{composable} if
$i(a) = t(b)$, in which case there exists an edge $c=ab$ of $X'$
such that $i(c) = i(b)$, $t(c) = t(a)$ and $a$, $b$ and $c$ form the
boundary of a $2$--simplex in $X'$.

A \emph{complex of groups} $G(X)=(G_\sigma, \psi_a, g_{a,b})$ over a
polyhedral complex $X$ is given by: \begin{enumerate} \item a group
$G_\sigma$ for each $\sigma \in V(X')$, called the \emph{local
group} at $\sigma$; \item a monomorphism $\psi_a\co
G_{i(a)}\rightarrow G_{t(a)}$ for each $a \in E(X')$; and
\item for each pair of composable edges $a$, $b$ in $X'$, an element $g_{a,b} \in
G_{t(a)}$, such that \[ \Ad(g_{a,b})\circ\psi_{ab} = \psi_a
\circ\psi_b
\] where $\Ad(g_{a,b})$ is conjugation by $g_{a,b}$ in $G_{t(a)}$,
and for each triple of composable edges $a,b,c$ the following
cocycle condition holds
\[\psi_a(g_{b,c})\,g_{a,bc} = g_{a,b}\,g_{ab,c}\] \end{enumerate} All complexes of
groups in this paper will be \emph{simple}, meaning that each
$g_{a,b}$ is trivial.

Next we define morphisms of complexes of groups. Let
$G(X)=(G_\sigma, \psi_a)$ and $H(Y)=(H_\tau,\psi_b)$ be simple
complexes of groups over polyhedral complexes $X$ and $Y$. Let $f\co
X'\to Y'$ be a simplicial map sending vertices to vertices and edges
to edges (such an $f$ is \emph{nondegenerate}). A \emph{morphism}
$\phi\co G(X) \to H(Y)$ over $f$ consists of:
\begin{enumerate}
\item a homomorphism $\phi_\sigma\co G_\sigma \to H_{f(\sigma)}$ for
each $\sigma \in V(X')$; and \item an element $\phi(a) \in
H_{t(f(a))}$ for each $a \in E(X')$, such that \[
\Ad(\phi(a))\circ\psi_{f(a)}\circ\phi_{i(a)} =
\phi_{t(a)}\circ\psi_a\] and for all pairs of composable edges
$(a,b)$ in $E(X')$,
\[ \phi(ab) = \phi(a) \,\psi(\phi(b)) \]
\end{enumerate}
We note that morphisms may also be defined over degenerate maps
$f\co X' \to Y'$. If $f$ is an isometry and each $\phi_\sigma$ an
isomorphism then $\phi$ is an \emph{isomorphism}.  A morphism
$\phi\co G(X) \to H(Y)$ is a \emph{covering} if:
\begin{enumerate}\item each $\phi_\sigma$ is injective;
and \item \label{i:covbijection} for each $\sigma \in V(X')$ and
$b \in E(Y')$ such that $t(b) = f(\sigma)$, the map
\[  \coprod_{\substack{a \in f^{-1}(b)\\ t(a)=\sigma}}
G_\sigma / \psi_a(G_{i(a)}) \to H_{f(\sigma)} / \psi_b(H_{i(b)})\]
induced by $g \mapsto \phi_\sigma(g)\phi(a)$ is a
bijection.\end{enumerate}

Let $G$ be a group acting without inversions on a polyhedral complex
$Y$.  This action induces a complex of groups, as follows. Let $X =
G \bs Y$ with $p\co Y \rightarrow X$  the natural projection. For
each $\sigma \in V(X')$, choose $\tilde\sigma \in V(Y')$ such that
$p(\tilde\sigma) = \sigma$. The local group $G_\sigma$ is the
stabiliser of $\tilde\sigma$ in $G$, and the $\psi_a$ and $g_{a,b}$
are defined using further choices. The resulting complex of groups
$G(X)$ is unique up to isomorphism. A complex of groups is
\emph{developable} if it is isomorphic to a complex of groups
associated, as just described, to an action.

Let $G(X)$ be a (simple) complex of groups, and choose a basepoint
$\sigma_0 \in V(X')$.  The \emph{fundamental group of the complex of
groups}
  $\pi_1(G(X),\sigma_0)$ is defined so that
if $X$ is simply connected, $\pi_1(G(X), \sigma_0)$ is isomorphic to
the direct limit of the family of groups $G_\sigma$ and
monomorphisms $\psi_a$.

If $G(X)$ is developable, then it has a \emph{universal cover}
$\widetilde{G(X)}$.  This is a simply-connected polyhedral complex,
equipped with an action of $\pi_1(G(X),\sigma_0)$, so that the
complex of groups induced by the action of the fundamental group on
the universal cover is isomorphic to $G(X)$.

We now describe a geometric condition for developability.  Let $X$
be a connected polyhedral complex and $\sigma \in V(X')$.  The
\emph{star} of $\sigma$, written $\St(\sigma)$, is the union of the
interiors of the simplices in $X'$ which meet $\sigma$.  If $G(X)$
is a complex of groups over $X$ then, even if $G(X)$ is not
developable, each $\sigma \in V(X')$ has a \emph{local development}.
That is, we may associate to $\sigma$ an action of $G_\sigma$ on the
star $\St(\tilde\sigma)$ of a vertex $\tilde\sigma$ in some
simplicial complex, such that $\St(\sigma)$ is the quotient of
$\St(\tilde\sigma)$ by the action of $G_\sigma$.  If $G(X)$ is
developable, then for each $\sigma \in V(X')$, the local development
of $\sigma$ is isomorphic to the star of each lift $\tilde\sigma$ of
$\sigma$ in the universal cover $\widetilde{G(X)}$.

The local development
$\St(\tilde\sigma)$ has a metric structure induced by that of the
polyhedral complex $X$.  We say that a complex of groups $G(X)$ is
\emph{nonpositively curved} if for all $\sigma \in V(X')$,
$\St(\tilde\sigma)$ has nonpositive curvature (that is,
$\St(\tilde\sigma)$ is locally CAT$(\kappa)$ for some $\kappa \leq
0$) in this induced metric. The importance of this condition is
given by:

\begin{theorem}[Haefliger \cite{bh1:ms}]
\label{t:nonpos}
A nonpositively curved complex of groups is developable.
\end{theorem}

Let $G(X)$ be a developable complex of groups over a polyhedral
complex $X$, with universal cover $Y$ and fundamental group $\G$. We
say that $G(X)$ is \emph{faithful} if the action of $\G$ on $Y$ is
faithful.  If $G(X)$ is faithful, then $\G$ may be regarded as a
subgroup of $\Aut(Y)$.  Moreover, $\Gamma$ is discrete if and only
if all local groups of $G(X)$ are finite, and a discrete $\G$ is a uniform
lattice if and only if $X$ is a finite polyhedral complex.

%%%%%%%%%%%%%%%%%%%%%%%%%%%%%%%%%%%%%%%%%%%%%%%%%%%%%%%%%%%%%%%%%%%%%
\subsubsection{Towers}\label{ss:towers}
%%%%%%%%%%%%%%%%%%%%%%%%%%%%%%%%%%%%%%%%%%%%%%%%%%%%%%%%%%%%%%%%%%%%%

Let $G(X)=(G_\sigma, \psi_a)$ and $H(X) = (H_\tau, \psi_b)$ be
simple complexes of groups over a complete, connected polyhedral
complex $X$. We say that $G(X)$ is a \emph{full complex of
subgroups} of $H(X)$ if there is a covering $\phi\co G(X) \to H(X)$
over the identity map $X' \to X'$ such that each $\phi(a) = 1$.  By
\cite[Corollary~3.16, Chapter~III.\cG]{bh1:ms} the covering $\phi$
induces an injective homomorphism of fundamental groups
\[\pi_1(G(X),\sigma_0) \to \pi_1(H(X),\sigma_0)\]

Suppose now that $Y$ is a polyhedral complex, and that we have a
sequence $(G_i(X))$ of complexes of groups over $X$, such that
\begin{enumerate} \item for each $i$, $G_i(X)$ is a full complex of subgroups of
$G_{i+1}(X)$ with respect to a covering $\phi_i$;
\item the image of each $(\phi_i)_\sigma$ is a proper subgroup; and
\item $Y$ is the universal cover of each
$G_i(X)$.\end{enumerate} The sequence $G_i(X)$ then induces an
infinite strictly ascending chain of fundamental groups
\[ \pi_1(G_1(X),\sigma_0) < \pi_1(G_2(X),\sigma_0) < \cdots <
\pi_1(G_i(X),\sigma_0) < \cdots
\] If each $G_i(X)$ is a faithful complex of finite groups, with finite covolume,
then each fundamental group in this chain is a lattice in $\Aut(Y)$.
In this way we may construct towers of lattices in $\Aut(Y)$.

%%%%%%%%%%%%%%%%%%%%%%%%%%%%%%%%%%%%%%%%%%%%%%%%%%%%%%%%%%%%%%%%%%%%
%%%%%%%%%%%%%%%%%%%%%%%%%%%%%%%%%%%%%%%%%%%%%%%%%%%%%%%%%%%%%%%%%%%
\section[The Functor Theorem]{The \fullref{thm:functor}}\label{s:extend}
%%%%%%%%%%%%%%%%%%%%%%%%%%%%%%%%%%%%%%%%%%%%%%%%%%%%%%%%%%%%%%%%%%%%%
%%%%%%%%%%%%%%%%%%%%%%%%%%%%%%%%%%%%%%%%%%%%%%%%%%%%%%%%%%%%%%%%%%%%%

This section contains our key technical result, the \fullref{thm:functor}
below. Let $\cG$ be the category of graphs of groups and morphisms,
as defined by Bass~\cite{b1:ctgg}. Let $\cC$ be the category of
complexes of groups and morphisms, as defined by
Haefliger~\cite{bh1:ms}. In \fullref{ss:GandC1} we make precise
the relationship between $\cG$ and the subcategory of $\cC$
consisting of complexes of groups over $1$--dimensional polyhedral
complexes. This is used in \fullref{ss:functorthm} where we
state and prove the \fullref{thm:functor}.

%%%%%%%%%%%%%%%%%%%%%%%%%%%%%%%%%%%%%%%%%%%%%%%%%%%%%%%%%%%%%%%%%%%%%%
%%%%%%%%%%%%%%%%%%%%%%%%%%%%%%%%%%%%%%%%%%%%%%%%%%%%%%%%%%%%%%%%%%%%%%
\subsection{Relationship between graphs of groups and $1$--dimensional
complexes of groups}\label{ss:GandC1}
%%%%%%%%%%%%%%%%%%%%%%%%%%%%%%%%%%%%%%%%%%%%%%%%%%%%%%%%%%%%%%%%%%%%%
%%%%%%%%%%%%%%%%%%%%%%%%%%%%%%%%%%%%%%%%%%%%%%%%%%%%%%%%%%%%%%%%%%%%%

Complexes of groups are generalisations of graphs of groups.  However,
Haefliger remarks~\cite[page~566]{bh1:ms} that it is unclear whether
the notion of morphism is the same.  Let $\cC_1$ be the subcategory
of $\cC$ consisting of complexes of groups over $1$--dimensional
polyhedral complexes (that is, simplicial graphs), and morphisms over
\emph{nondegenerate} polyhedral maps. The following result specifies
the relationship between the category $\cG$ of graphs of groups and the
category $\cC_1$.

\begin{proposition}\label{p:gtoc1} There is a functor $F\co\cG \to \cC_1$ which is
a bijection on objects, takes faithful graphs of groups to faithful
complexes of groups, surjects onto the set of morphisms of $\cC_1$
and takes coverings to coverings. However, there is no functor from
$\cG$ to $\cC_1$ which is injective on the set of morphisms of
$\cG$.  Hence these categories are not isomorphic.
\end{proposition}

\begin{proof} We first define the functor $F$ on objects.
Let $\mathbb{A} = (A,\cA)$ be a graph of groups and let $|A|$ be the
geometric realisation of the graph $A$.  We construct a complex of
groups $F(\mathbb{A})$ over $|A|$.  The local groups at the vertices
of $|A|$ are the vertex groups of $\mathbb{A}$.  For each $e \in
E(A)$, let $\sigma_e=\sigma_{\overline{e}}$ be the vertex of the
barycentric subdivision $|A|'$ at the midpoint of $e$. Then the
local group at $\sigma_e$ in the complex of groups is $\cA_e =
\cA_{\overline{e}}$. Each monomorphism $\alpha_e\co \cA_e \to
\cA_{i(e)}$ of the graph of groups induces the same monomorphism of
corresponding local groups in the complex of groups $F(\mathbb{A})$.

It is clear that $F$ is a bijection from the
set of objects of $\cG$ to the set of objects of $\cC_1$. From the definitions
of the actions on the respective universal covers, one sees also that faithful
graphs of groups are mapped to faithful complexes of groups.

To define the functor $F$ on morphisms, let $\phi\co\mathbb{A} \to
\mathbb{B}$ be a morphism of graphs of groups over a graph morphism
$f\co A \to B$.  By abuse of notation $f$ induces a nondegenerate
polyhedral map $f\co|A|' \to |B|'$.  We define a morphism of
complexes of groups $F(\phi)\co F(\mathbb{A}) \to F(\mathbb{B})$
over $f$.  The homomorphisms on local groups are the same as for the
morphism $\phi$. Let $a$ be an edge of $|A|'$. Then the monomorphism
$\psi_a$ in the complex of groups $F(\mathbb{A})$ is the same as
some $\alpha_e\co\cA_e \to \cA_{i(e)}$ in the graph of groups
$\mathbb{A}$. Put $\phi(a) = \delta_e$, and we obtain the required
\[ \phi_{t(a)}\circ \psi_a = \Ad(\phi(a)) \circ \psi_{f(a)}\circ
\phi_{i(a)}\] Hence $F(\phi)$ is a morphism of complexes of groups.
The map $F$ respects composition of morphisms, and is thus a
functor.

Given a morphism of $\cC_1$, one can construct a morphism of the
corresponding objects in $\cG$ by setting all $\gamma_v = 1$, and
all $\gamma_e = \phi(a)$, where the monomorphism $\alpha_e$ is the
same as the monomorphism $\psi_a$. Hence $F$ surjects onto the set
of morphisms of $\cC_1$. By checking the definition of covering in
both categories, we find that $F$ sends coverings to coverings, and
surjects onto the set of coverings of $\cC_1$.

However, since many choices of $\gamma_v$ and $\gamma_e$ could lead to the
same collection of $\delta_e$, there is no functor from $\cG$ to
$\cC_1$ which is injective on the set of morphisms of $\cG$.
\end{proof}

%%%%%%%%%%%%%%%%%%%%%%%%%%%%%%%%%%%%%%%%%%%%%%%%%%%%%%%%%%%%%%%%
\subsection[The Functor Theorem]{The \fullref{thm:functor}}\label{ss:functorthm}
%%%%%%%%%%%%%%%%%%%%%%%%%%%%%%%%%%%%%%%%%%%%%%%%%%%%%%%%%%%%%%%%

Let $T$ be a tree and $X$ a right-angled building.  We write
$\cG(T)$ for the category of graphs of groups with universal
covering tree $T$ and $\cC(X)$ for the category of developable
complexes of groups with universal cover $X$.  We now construct a
functor from $\cG(T)$ to $\cC(X)$, for certain $T$ and $X$.

This functor will exist only when the building $X$ is ``sufficiently
symmetric", as defined below. However, the proof of the
\fullref{thm:functor} will show that even if $X$ is not sufficiently symmetric, we
can construct a functor taking faithful graphs of groups to faithful
complexes of groups, and coverings to coverings, if we restrict to
graphs of groups whose vertices are $2$--colourable and to coverings
which preserve this colouring. As shown in \fullref{s:proofmain}
below, the proof of the \fullref{thm:main} requires only this restricted
functor.

To explain when a right-angled building $X$ is sufficiently
symmetric, let $(W,I)$ be the right-angled Coxeter system, $P'$ the
chamber and $\{q_i\}$ the parameters associated to $X$.  Suppose
$m_{i_1,i_2}=\infty$, and consider the following symmetry
conditions:
\begin{enumerate}
\item \label{i:isomPP} there exists a bijection $g$ on the set $I$
such that $m_{i,j}=m_{g(i),g(j)}$ for all $i,j \in I$, and
$g(i_1)=i_2$.  This may be thought of as an isometry of $P'$ which
takes the $i_1$--face to the $i_2$--face.
\item \label{i:isomface} there exists a bijection\[ h \co \{i \in I
\,|\, m_{i_1,i} < \infty \} \to \{i \in I \,|\, m_{i_2,i} < \infty
\}
\] such that $m_{i,j}=m_{h(i),h(j)}$ for all $i,j$ in the domain,
$h(i_1)=i_2$, and $q_{i}=q_{h(i)}$ for all $i$ in the domain.  This
may be thought of as an isometry from the simplicial neighbourhood
in $P'$ of the $i_1$--face to the simplicial neighbourhood in $P'$
of the $i_2$--face, which preserves the cardinalities $q_i$.
\end{enumerate}

Conditions~\eqref{i:isomPP} and~\eqref{i:isomface} are satisfied if,
for example, $P'$ is the barycentric subdivision of a regular
right-angled hyperbolic polyhedron $P$, and all $q_i$ are equal. In
dimension 2, $P$ may be a $k$--gon for any $k \geq 5$. In each of
dimensions 3 and 4, there is only one such $P$, the dodecahedron and
the $120$--cell respectively (see, for example, Vinberg and
Shvartsman~\cite{vs1:g}).

Conditions~\eqref{i:isomPP} and~\eqref{i:isomface} may also be
satisfied for $P$ which is not regular.  For example, suppose $P$ is
$3$--dimensional.  We claim that there exist $i_1, i_2 \in I$ such
that $m_{i_1,i_2}=\infty$ and the $i_1$-- and $i_2$--faces of $P$
are isometric.  Since the $2$--dimensional faces of $P$ are regular,
it suffices to exhibit a pair of nonadjacent $2$--dimensional faces
with, say, $5$ sides.  Let $a_2$ be the number of $2$--dimensional
faces of $P$. For each $2$--dimensional face $F$, let $a_1(F)$ be
the number of sides of $F$ and let $\ex(F) = a_1(F) - 5 \geq 0$.
Elementary calculations yield
\[\label{e:pv} a_2 = 12 + \!\sum_{\dim(F)=2} \ex(F)\]
(see Potyagailo and Vinberg~\cite[page 70]{pv1:rrg}).
This implies that  there are at least $12$ faces $F$ with $5$ sides.
Hence we can find two nonadjacent $2$--dimensional faces of $P$ with
$5$ sides. This kind of argument cannot be extended to dimension
$4$, since the $3$--dimensional faces of $P$ will not in general be
regular.

We now state and prove the \fullref{thm:functor}.

\begin{functortheorem}
\label{thm:functor}
Let $X$ be a right-angled building of
type $(W,I)$ and parameters $\{q_i\}$.
  For each $i_1, i_2 \in I$ such that $m_{i_1,i_2} =
\infty$, let $T$ be the $(q_{i_1},q_{i_2})$--biregular tree.
Suppose condition~\eqref{i:isomPP} above holds, and that if $q_{i_1}
= q_{i_2}$ then condition~\eqref{i:isomface} above holds, with $g$
an extension of $h$.  Then there is a functor $F$ from $\cG(T)$ to
$\cC(X)$, which takes faithful graphs of groups to faithful
complexes of groups, and coverings to coverings.
\end{functortheorem}

\begin{proof} By \fullref{p:gtoc1}, it suffices to construct
a functor $F$ from the image of $\cG(T)$ in $\cC_1$ to $\cC(X)$.  We
first define $F$ on objects.  Let $G(Y)$ be an object of $\cC_1$
which is in the image of $\cG(T)$.

For each edge $e$ of the simplicial graph $Y$, let $P'_e$ be a copy
of $P'$ (recall that $P'$ is the simplicial cone with vertex $x_0$
on the barycentric subdivision of the finite nerve of $(W,I)$).
Identify the midpoint of $e$ with the vertex $x_0$ of $P'_e$.

Suppose first that the vertices of the graph $Y$ may be
$2$--coloured with the types $i_1$ and $i_2$.  If $q_{i_1}= q_{i_2}$
choose such a $2$--colouring, and if $q_{i_1}\not = q_{i_2}$, then
assign the colours $i_1$ and $i_2$ according to the valences of the
vertices in the universal covering tree $T$. Then for $j = 1,2$
identify the vertex of $e$ of type $i_j$ with the $i_j$--vertex of
$P'_e$.

Now suppose the vertices of the graph $Y$ are not $2$--colourable.
If the edge $e$ does not form a loop in $Y$, then identify one
vertex of $e$ with the $i_1$--vertex of $P'_e$, and the other vertex
with the $i_2$--vertex. If the edge $e$ does form a loop in $Y$ then
use the isometry $h$ of condition~\eqref{i:isomface} above to form
$P'_e/h$. That is, glue together the $i_1$-- and $i_2$--faces of
$P'_e$ using $h$.  Then identify the vertex of $Y$ to which both
ends of $e$ are joined with the image of the $i_1$-- and
$i_2$--vertices in $P'_e/h$.

Now, glue together, either by preserving type on $i_1$-- and $i_2$--faces,
or if the vertices of $Y$ are not $2$--colourable, by using the isometry
$h$, the faces of the various $P'_e$ and $P'_e/h$ whose centres
correspond to the same vertex of $Y$. Let $F(Y)$ be the resulting
polyhedral complex. By construction, the barycentric subdivision $Y'$
embeds in $F(Y)$. Note also that each vertex of $F(Y)$ has either the
type of a unique $J \in \mathcal{S}^f$, or the two types $J$ and $h(J)$,
where $i_1 \in J$ and $J \in \mathcal{S}^f$.

We now construct a complex of groups $F(G(Y))$ over $F(Y)$.  First,
fix the local groups and monomorphisms induced by the embedding of
$Y'$ in $F(Y)$. For each $i \in I$, let $G_i = \Z / q_i \Z$.  For
each edge $e$ of $Y$, let $G_e$ be the local group at the midpoint
of $e$ in the complex of groups $G(Y)$.

Consider $J \in \mathcal{S}^f$ which does not contain $i_1$ or $i_2$. The
local group at the vertex with type $J$ of $P'_e$ or $P'_e/h$ is
then
\[G_e \times \prod_{j \in J} G_j\]  The monomorphisms between
such local groups are natural inclusions.

Now consider $J \in \mathcal{S}^f$ which contains one of $i_1$ or $i_2$ (since
$m_{i_1,i_2} = \infty$, $J$ cannot contain both $i_1$ and $i_2$).
Without loss of generality suppose $i_1 \in J$, and let $F_e$ be
either the $i_1$--face of $P'_e$, or the glued face of $P'_e/h$. The
vertex of type $J$ in $P'_e$ or $P'_e/h$ is then contained in $F_e$.
Let $v$ be the vertex of $Y$ which is identified with the centre of
$F_e$, and let $G_v$ be the local group at $v$ in the complex of
groups $G(Y)$. The local group at the vertex of type $J$ of $P'_e$
or $P'_e/h$ is then \[ \label{e:localgp}G_v \times
\prod_{\substack{j \in J
\\ j \not = i_1}} G_j\] This is well-defined, since if the face of
type $J$ also has type $h(J)$, by condition~\eqref{i:isomface} we
have $q_j = q_{h(j)}$ for each $j \in J$, hence $G_j = G_{h(j)} =
\Z/ q_j \Z$. The monomorphism from $G_v$ to this local group is the
natural inclusion onto the first factor. For each $J' \subset J$
with $i_1 \in J'$, the monomorphism
\[G_v \times \prod_{\substack{j \in J' \\ j \not = i_1}} G_j
\longrightarrow G_v \times \prod_{\substack{j \in J \\ j \not =
i_1}} G_j\] is the natural inclusion.  For each $J' \subset J$ with
$i_1 \not \in J'$, the monomorphism
\[G_e \times \prod_{j \in J'} G_j \longrightarrow G_v
\times \prod_{\substack{j \in J \\ j \not = i_1}} G_j\] is a
monomorphism $G_e \to G_v$ from the complex of groups $G(Y)$ on the
first factor, and natural inclusion on the other factors.

We now show that the complex of groups $F(G(Y))$ is developable,
with universal cover $X$. Let $\tau$ be a vertex of $F(Y)$ of type
$J$, with $|J|=n=\dim(X)$. One verifies that the link in the local
development at $\tau$ is the join of $n$ sets of points of
cardinalities respectively $q_j$, for $j \in J$.  Thus, as $P'$ is
right-angled, by Gromov's Link Condition (see Bridson and
Haefliger~\cite[Theorem~5.2, Chapter~II]{bh1:ms}) each local development
has nonpositive curvature. It follows by \fullref{t:nonpos} above that
the complex of groups $F(G(Y))$ is developable. The universal cover is
a building (see Gaboriau and Paulin~\cite[Section~3.3]{gp1:ih}), and by
\fullref{t:rabeu} above, the universal cover is the unique right-angled
building $X$ of type $(W,I)$ and parameters $\{q_i\}$.

The universal covering tree $T$ naturally embeds in $X$, and by
construction, if the fundamental group of $G(Y)$ acts faithfully on
$T$, then the fundamental group of $F(G(Y))$ acts faithfully on $T$,
and hence on $X$.  Thus $F$ maps faithful graphs of groups to
faithful complexes of groups.

We now define the functor $F$ on morphisms.  Let $G(Y)$ and $H(Z)$
be objects of $\cC_1$ in the image of $\cG(T)$, and let $\phi\co
G(Y) \to H(Z)$ be a morphism over a nondegenerate map $f\co Y' \to
Z'$. Let $F(Y)$ and $F(Z)$ be the polyhedral complexes constructed
from the simplicial graphs $Y$ and $Z$ as above.

If the vertices of $Y$ and $Z$ are $2$--colourable by the types
$i_1$ and $i_2$, and $f$ preserves these types, then $f$ may be
extended to a polyhedral map $F(f)\co F(Y) \to F(Z)$ by preserving
type on each copy of $P'$. Otherwise, we use
condition~\eqref{i:isomPP} to construct $F(f)\co F(Y) \to F(Z)$.  On
each copy of $P'$ or $P'/h$, $F(f)$ maps the vertex of type $J$ to
the vertex of type $g(J)$. This is well-defined since $h$ is the
restriction of $g$.

We now construct a morphism $F(\phi)$ of complexes of groups over
$F(f)$.  If $\tau$ is a vertex of $F(Y)$ then the local group at
$\tau$ is a direct product
\[G_\sigma \times \prod G_j \]
where $\sigma$ is a vertex of $Y'$.  The homomorphism of local
groups
 \[
G_\sigma \times \prod G_j \longrightarrow H_{f(\sigma)} \times \prod
G_j\] is $\phi_\sigma$ on the first factor, and the identity on
other factors.  Let $b$ be an edge of $F(Y)$.  If $\psi_b$, the
monomorphism along the edge $b$ in the complex of groups $F(G(Y))$,
has as its first factor a monomorphism $\psi_{a}$ from the complex
of groups $G(Y)$, put $F(\phi)(b)=\phi(a)$.  Otherwise, put
$F(\phi)(b) = 1$.

The following claims all follow from definitions: $F(\phi)$ is a
morphism of complexes of groups, $F(\phi)$ respects composition, and
$F(\phi)$ takes coverings to coverings. This completes the proof of
the \fullref{thm:functor}.
\end{proof}

%%%%%%%%%%%%%%%%%%%%%%%%%%%%%%%%%%%%%%%%%%%%%%%%%%%%%
%%%%%%%%%%%%%%%%%%%%%%%%%%%%%%%%%%%%%%%%%%%%%%%%%%%%%%%%%%%%%%%%%%%
\section[Proof of the Main Theorem]{Proof of the \fullref{thm:main}}
\label{s:proofmain}
%%%%%%%%%%%%%%%%%%%%%%%%%%%%%%%%%%%%%%%%%%%%%%%%%%%%%%%%%%%%%%%%%%%
%%%%%%%%%%%%%%%%%%%%%%%%%%%%%%%%%%%%%%%%%%%%%%%%%%%%%%%%%%%%%%%%%%%

We conclude by proving the \fullref{thm:main}, stated in the Introduction.
Most parts of the theorem are proved using the \fullref{thm:functor}.
Throughout this section, $X$ is as in the statement of the
\fullref{thm:main}, and $G = \Aut(X)$.

For covolume results, we will use \fullref{c:covolcalc} below,
which follows from the proof of the \fullref{thm:functor}.  Recall that
for a polyhedral complex $Y$, we write $Y_n$ for the set of
$n$--dimensional cells of $Y$.

\begin{corollary}\label{c:covolcalc}
Suppose $G(Y)$ is a faithful complex of finite groups over an
$n$--dimensional polyhedral complex $Y$, with universal cover $X$.
Let $\Gamma$ be the fundamental group of the complex of groups
$G(Y)$. For each cell $\sigma \in Y_n$, let $\Gamma_\sigma$ be the
local group at the barycentre of that cell. Then if $G(Y) =
F(\mathbb{A})$, where $F$ is the functor defined in the
\fullref{thm:functor}, and $\mathbb{A}=(A,\cA)$ is a graph of groups,
\[\mu(\G\bs G) =   \sum_{\sigma \in Y_n}
\frac{1}{|\Gamma_\sigma|}  =  \sum_{e \in E(A)} \frac{1}{|\cA_e|}
\]
That is, $\mu(\G \bs G)$ equals the covolume of the graph of groups
$\mathbb{A}$.
\end{corollary}

%%%%%%%%%%%%%%%%%%%%%%%%%%%%%%%%%%%%%%%%%%%%%%%%%%%%%%%%%%%%%%%%%%%%%
\subsection{Uniform covolumes}\label{s:unifcovol}
%%%%%%%%%%%%%%%%%%%%%%%%%%%%%%%%%%%%%%%%%%%%%%%%%%%%%%%%%%%%%%%%%%%%%%

%%%%%%%%%%%%%%%%%%%%%%%%%%%%%%%%%%%%%%%%%%%%%%%%%%%%%%%%%%%%%%%%%%%%%%%%
\subsubsection[Proof of~\ref{i:unifcovol}]{Proof of~\eqref{i:unifcovol}}
%%%%%%%%%%%%%%%%%%%%%%%%%%%%%%%%%%%%%%%%%%%%%%%%%%%%%%%%%%%%%%%%%%%%%%%%

We first show that the covolume of every uniform lattice belongs to
the set $\cV_u(G)$ of rationals given in~\eqref{i:unifcovol}.  Then
we show that any $v \in \cV_u(G)$ is the covolume of some uniform
lattice.

Since $X$ is a building, the group $G=\Aut(X)$ has a finite index
normal subgroup $\Aut_0(X)$, the group of type-preserving
automorphisms, which acts without inversions.  Thus, any uniform
lattice $\G < G$ has a finite index subgroup $\G \cap \Aut_0(X)$
which acts without inversions. Now, the set $\cV_u(G)$ is closed
under multiplication by positive integers. Hence, to show that the
covolume of every uniform lattice belongs to $\cV_u(G)$, it suffices
to consider only uniform lattices which act without inversions.

Let $\G$ be a uniform lattice in $G$ which acts without inversions,
and let
\[ \mu(\G \bs G) = \Vol(\G \quot X_n) = \frac{a}{b}\]
To show that the prime divisors of $b$ are strictly less than
$\max_i \{q_i\}$, we use the following.  Let $x$ be a vertex of $X$,
of type some $J \in \cS^f$ with $|J|=n=\dim(X)$. Let $S_k$ be the
symmetric group on $k$ letters. The group of type-preserving
automorphisms of the link of $x$ in $X$ is then\[\prod_{j \in J}
S_{q_j}\] and its subgroup which fixes an $(n-1)$--dimensional cell
of the link pointwise is
\[\Fix(n-1,J):= \prod_{j \in J} S_{q_j-1}.\]
A result in a forthcoming paper~\cite[Theorem~3]{at1:cul} then
implies the restriction on the prime divisors of $b$; we sketch the
argument.  It suffices to show that for any $\sigma \in X_n$, the
prime divisors of $|\G_\sigma|$ are strictly less than $\max_i
\{q_i\}$. Since $\G_\sigma$ is a finite group, for some $m > 0$ it
injects into $H_m:=\Aut(B(m,\sigma))$, the group of automorphisms of
the combinatorial ball of radius $m$ in $X$ centred at $\sigma$.  It
now suffices to bound the prime divisors of $|H_m|$.  For this, we
use induction on $m$, and the fact that the prime divisors of the
group $\Fix(n-1,J)$ above are strictly less than $\max_i\{q_i\}$. We
conclude that the covolume of any uniform lattice in $G$ belongs to
the set $\cV_u(G)$.

We now show that every $v \in \cV_u(G)$ is the covolume of some
uniform lattice. By \fullref{c:covolcalc} above, and the proof
of the \fullref{thm:functor}, it suffices to show that, for any rational
$\frac{a}{b}$ with the prime divisors of $b$ strictly less than
$\max_i \{ q_i \}$, there is a faithful graph of finite groups
$\mathbb{A}=(A,\cA)$ over a finite graph $A$ such that
\begin{enumerate}
\item the universal covering tree of $\mathbb{A}$ is
$(q_{i_1},q_{i_2})$--biregular, where $m_{i_1,i_2}=\infty$,
\item the vertices of the graph $A$ may be $2$--coloured, and
\item the covolume of $\mathbb{A}$ equals $\frac{a}{b}$.
\end{enumerate}

Suppose first that all $q_i$ are equal, and let $T_q$ be the
$q$--regular tree.

If all $q_i = 2$, we construct a graph of groups with universal
covering tree $T_2$ and edge covolume $\frac{a}{1} =a $, for any
integer $a \geq 1$. Let $A$ be the graph which is a ray of $a$
edges.  The vertices of $A$ are then $2$--colourable.  Let
$\mathbb{A}=(A,\cA)$ be the graph of groups over $A$ with each edge
group trivial, the two vertex groups at either end of the ray both
$\Z / 2\Z$, and all other vertex groups trivial. Since the edge
groups are trivial this graph of groups is faithful.  Its universal
covering tree is $T_2$ and it has covolume $a$.

\begin{figure}
\begin{center}
\includegraphics[width=6cm]{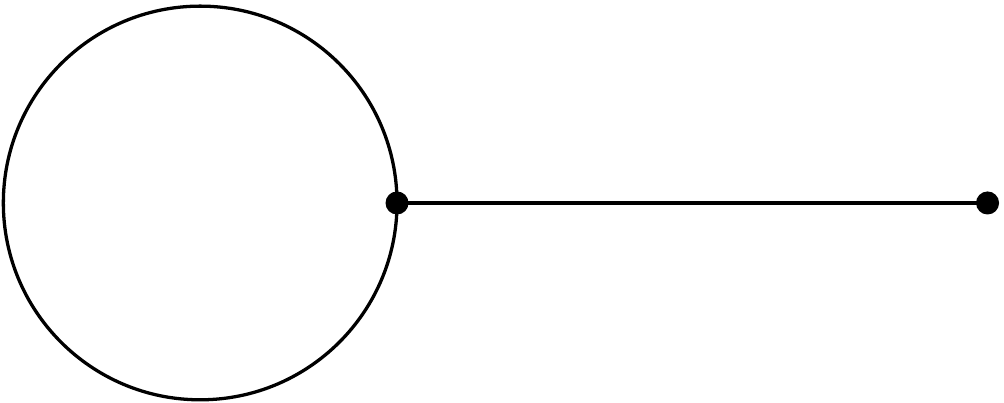}
%\begin{pspicture}(1,0.5)(8,4)
%\psline(4,2)(7,2)\psarc(3,2){1}{0}{360}%\psarc(3,2){1}{0}{181}
%\psdot*(4,2)\psdot*(7,2)
%\end{pspicture}
\end{center}
\caption{Rosenberg's graph for the case $q=3$}\label{f:graph}
\end{figure}

If all $q_i=q$ are equal and $q\geq 3$, we require a graph of groups
with universal covering tree $T_q$ and covolume $\frac{a}{b}$, where
the prime divisors of $b$ are strictly less than $q$. This is
exactly the set of covolumes of uniform lattices in $\Aut(T_q)$
(see Rosenberg~\cite[Proposition~9.1.2]{ros1:tctl}).  For $q \geq
4$, the vertices of each of the graphs in Rosenberg's proof are
$2$--colourable. In the case $q = 3$, though, some of the graphs of
groups are over the graph shown in \fullref{f:graph}.  The vertices of
this graph cannot be $2$--coloured.

A key idea of Rosenberg's proof, which is also used below to
prove~\eqref{i:countunif}, is that of finite-sheeted covers. In each
case, the graph $A$, considered as a topological space, is not
simply connected.  The fundamental group $\pi_1(A)$ is then a free
group on finitely many generators, so has subgroups of arbitrary
finite index.  Hence $A$ has a finite-sheeted topological cover for
any finite number of sheets.  Thus, if a uniform lattice constructed
over $A$ has covolume $\frac{a}{b}$, we may construct a uniform
lattice with covolume a positive integer multiple of $\frac{a}{b}$,
over a finite-sheeted cover of $A$.

In particular, the $2$--sheeted cover of the graph in
\fullref{f:graph} has 4 vertices, which may be $2$--coloured. As
explained in \cite[Example~4, Chapter~4]{ros1:tctl}, the edge groups
of this graph of groups have order a power of $2$, and this power
may be chosen arbitrarily. Thus, by taking a $2$--sheeted cover and
then doubling the order of the edge groups, we obtain a graph of
groups over a $2$--colourable graph with unchanged covolume. This
completes the proof of~\eqref{i:unifcovol} in the case where all
$q_i$ are equal.

Now suppose the $q_i$ are not all equal, and let
$q_{i_1}=\max_i\{q_i\}$.  By assumption, there is a $q_{i_2}$ with
$m_{i_1,i_2}=\infty$.  If $q_{i_1}=q_{i_2}$, then use the same
graphs of groups as in the case where all $q_i$ are equal.
Otherwise, to simplify notation, put $p = q_{i_1}$ and $q =
q_{i_2}$, so that $p > q \geq 2$, and let $T_{p,q}$ be the
$(p,q)$--biregular tree.

Rosenberg \cite[Theorem~9.2.1]{ros1:tctl} shows that the set of
covolumes of uniform lattices in $\Aut(T_{p,q})$ is the set of
rationals $\frac{a}{b}$ such that the prime divisors of $b$ are
strictly less than $p$.  As with regular trees, for each such
$\frac{a}{b}$, there is a faithful graph of finite groups, over a
finite graph, with covolume $\frac{a}{b}$. Since each graph is
covered by $T_{p,q}$, and $p \not = q$, the vertices of each graph
are $2$--colourable.  This completes the proof
of~\eqref{i:unifcovol}.

%%%%%%%%%%%%%%%%%%%%%%%%%%%%%%%%%%%%%%%%%%%%%%%%%%%%%%%%%%%%%%%%%%%%%%%
\subsubsection[Proof of~\ref{i:countunif}]{Proof of~\eqref{i:countunif}}
%%%%%%%%%%%%%%%%%%%%%%%%%%%%%%%%%%%%%%%%%%%%%%%%%%%%%%%%%%%%%%%%%%%%%%%

By the proof of the \fullref{thm:functor} and \fullref{c:covolcalc}, it suffices to prove
\fullref{p:treecount} below.

\begin{proposition}\label{p:treecount} Let $T$ be a regular or biregular tree,
other than the $2$--regular tree.  Let $v > 0$ be the covolume of a
uniform lattice in $\Aut(T)$.  Then there is a countably infinite
number of nonconjugate uniform lattices in $\Aut(T)$ which have
covolume $v$.
\end{proposition}

\begin{proof}  Since there are only countably many finite
graphs of finite groups, there are at most countably many uniform lattices in
$\Aut(T)$.

Assume that $T$ is the $m$--regular tree, for $m \geq 3$ (the proof for
biregular trees is similar). Rosenberg \cite[Proposition~9.1.2]{ros1:tctl}
constructs a faithful graph of
finite groups $\mathbb{A}=(A,\cA)$ over a finite graph $A$ such that
\begin{enumerate}
\item the universal covering tree of $\mathbb{A}$
is $T$, \item the  fundamental group $\G$ of $\mathbb{A}$ has
covolume $v$, \item the vertices of the graph $A$ may be
$2$--coloured (in the case where $T$ is $3$--regular, use a
$2$--sheeted cover),
\item \label{i:nsc} the graph $A$ is not simply connected, and \item
\label{i:multedge} the orders of the edge and vertex groups $\cA_e$
and $\cA_v$ may be multiplied by powers of primes strictly less than
$m$, and the new graph of groups is faithful.
\end{enumerate}

Let $n$ be a positive integer whose prime divisors are strictly less
than $m$. By~\eqref{i:nsc}, we may obtain a uniform lattice in
$\Aut(T)$ of covolume $nv$. Applying~\eqref{i:multedge}, we may
obtain a uniform lattice in $\Aut(T)$ of covolume $v/n$.  By
carrying out these two processes at the same time, we obtain a new
lattice $\G'$ of covolume $v$, which is nonconjugate to the original
lattice $\G$, since the graph $\G' \bs T$ is not the same as $A$.
This may be done for countably many values of $n$. Hence we obtain
countably many nonconjugate uniform lattices of covolume $v$.
\end{proof}

%%%%%%%%%%%%%%%%%%%%%%%%%%%%%%%%%%%%%%%%%%%%%%%%%%%%%%%%%%%%%%%%%%%%%%%%%%%
\subsection{Nonuniform covolumes}
%%%%%%%%%%%%%%%%%%%%%%%%%%%%%%%%%%%%%%%%%%%%%%%%%%%%%%%%%%%%%%%%%%%%%%%%%%%

%%%%%%%%%%%%%%%%%%%%%%%%%%%%%%%%%%%%%%%%%%%%%%%%%%%%%%%%%%%%%%%%%%%%%%%%%%%%
\subsubsection[Proof of~\ref{i:nonunifcovol}]{Proof
of~\eqref{i:nonunifcovol}}
%%%%%%%%%%%%%%%%%%%%%%%%%%%%%%%%%%%%%%%%%%%%%%%%%%%%%%%%%%%%%%%%%%%%%%%%%%%%

Let $T$ be a regular or biregular tree, other than the $2$--regular
tree, and let $v > 0$.   By \fullref{c:covolcalc} and the
proof of the \fullref{thm:functor}, it suffices to construct a nonuniform
lattice $\G$ in $\Aut(T)$ of covolume $v$, such that the vertices of
the graph $\G \bs T$ may be $2$--coloured.

Bass--Lubotzky showed that for $m \geq 3$ and every $v > 0$
there is a nonuniform lattice $\G$ acting on the $m$--regular
tree $T_m$, such that the covolume of $\G$ is $v$ (see Bass and
Lubotzky \cite[Theorem~4.3]{bl1:tl}).  Moreover, the graph $\G \bs
T_m$ is a ray, hence its vertices are $2$--colourable.  Rosenberg
\cite[Theorem~8.2.2]{ros1:tctl} showed the analogous result for biregular
trees, without the conclusion that the quotient graph is a ray. Since
the vertices of any graph covered by a biregular tree $T_{m,n}$,
with $m \not = n$, are $2$--colourable, this completes the proof
of~\eqref{i:nonunifcovol} of the \fullref{thm:main}.

%%%%%%%%%%%%%%%%%%%%%%%%%%%%%%%%%%%%%%%%%%%%%%%%%%%%%%%%%%%%%%%%%%%%%%%%%%%%
\subsubsection[Proof of~\ref{i:nonunifcount}]{Proof
of~\eqref{i:nonunifcount}}
%%%%%%%%%%%%%%%%%%%%%%%%%%%%%%%%%%%%%%%%%%%%%%%%%%%%%%%%%%%%%%%%%%%%%%%%%%%%

Let $T$ be the $(m,n)$--biregular tree, with $m,n \geq 3$.  Farb and
Hruska~\cite{fh1:ci} constructed commensurability invariants for
nonuniform lattices acting on $T$.  One such invariant is growth type.

As in \cite[Section~3.2]{fh1:ci}, let $f, g\co\N \to \N$ be any
two functions. We say that $f \preceq g$ if for some $\alpha, \beta
\in \N$,
\[ f(k) \leq \alpha\, g(k + \beta) \]
and we say that $f$ and $g$ are \emph{equivalent} if $f \preceq g$
and $g \preceq f$.

Let $Y$ be a locally finite graph with basepoint $\ast$, and let
$g\co\N \to \N$ be such that $g(k)$ is the number of vertices in the
combinatorial ball of radius $k$ centred at $\ast$.  The
\emph{growth type} of $(Y,\ast)$ is defined to be the equivalence
class of the function $g$.  The growth type is independent of the
choice of basepoint.  If $\G$ and $\G'$ are commensurable in
$\Aut(T)$, then using finite covers one sees that the graphs $\G \bs
T$ and $\G' \bs T$ have the same growth type
\cite[Proposition~3.4]{fh1:ci}.

The growth type of a lattice $\G$ in $\Aut(X)$ may be similarly
defined.  The ball of radius 1 around a basepoint $\ast$ in $\G \bs
X$ consists of all images of chambers which meet $\ast$, and by
induction the ball of radius $k$ consists of all images of chambers
which meet the ball of radius $k - 1$. This growth type is also a
commensurability invariant, by a similar proof to that for tree
lattices.

Using growth type, for any $v > 0$, there are uncountably many
commensurability classes of nonuniform lattices with covolume $v$ in
$\Aut(T)$ \cite[Theorem~5.2]{fh1:ci}.  These lattices are
constructed using graphs of groups over particular trees.  We now
show that $F$ in the \fullref{thm:functor}, at least when applied to a
tree, preserves growth type:

\begin{lemma}\label{l:growth} Let $Y$ be an (infinite) tree and let
$F$ be as in the \fullref{thm:functor}.  Then the growth type of $Y$ is
the same as the growth type of the polyhedral complex $F(Y)$.
\end{lemma}

\begin{proof} Choose a basepoint $\ast$ of $Y$ and let $g(k)$ be the number of vertices in the
combinatorial ball of radius $k$ centred at $\ast$.  Let $G(k)$ be
the number of vertices in the combinatorial ball of radius $k$
centred at the image of $\ast$ in $F(Y)$.  We claim that $g$ is
equivalent to $G$.

Let $l$ be the number of vertices of $P'$, and for $j = 1,2$ let
$l_j$ be the number of vertices of the $i_j$--face of $P'$. Without
loss of generality we may assume that $l_1 \geq l_2$ and that $\ast$
has type $i_1$.  Put $C = l_1 - l_2 \geq 0$.  Since $Y$ is a tree,
we then compute
\begin{equation}\label{e:Gg} G(k) = (l - l_j)\,g(k) - l + 2l_1 +
C\!\sum_{m=1}^{k-1}(-1)^m g(m) \end{equation} where $j=1$ if $k$ is
odd and $j = 2$ if $k$ is even.

To show $G \preceq g$, by induction on $k$ we have
\begin{equation}\label{e:bound} \sum_{m=1}^{k-1}(-1)^m g(m) \,\leq \,g(k)\end{equation}
Choose $\alpha$ so that $l - l_j + C \leq \alpha$ for $j = 1,2$.
Then choose $\beta$ so that $-l + 2l_1 \leq \alpha\beta$. As $Y$ is
an infinite tree, $\beta \leq g(k + \beta) - g(k)$ for all $k$.
Hence
\[ G(k) \leq \alpha g(k) - l + 2l_1  \leq \alpha g(k + \beta) \]
This proves $G \preceq g$.

To show $g \preceq G$, by~\eqref{e:Gg} we have \[ g(k) = \frac{1}{l
- l_j} \, G(k) + \frac{l - 2l_1}{l - l_j} - \frac{C}{l - l_j}\!
\sum_{m=1}^{k-1}(-1)^m g(m) \] If $k$ is odd, then by induction
\[  - \sum_{m=1}^{k-1}(-1)^m g(m) = g(1) - g(2) + \cdots -g(k-1) \leq g(1) + 0 \]
hence \begin{equation}\label{e:kodd} g(k) \leq \frac{1}{l - l_1} \,
G(k) + \frac{l - 2l_1}{l - l_1} + \frac{C}{l - l_1}\,g(1)
\end{equation}
 If $k$ is even then a similar inequality to~\eqref{e:bound}
yields \begin{equation}\label{e:keven} g(k) \leq  \frac{1}{l - l_2}
\, G(k) + \frac{l - 2l_1}{l - l_2} + \frac{C}{l - l_2}[g(1) + g(k)]
\end{equation}
Since $C= l_1 - l_2$, we have $1 - \frac{C}{l - l_2} = \frac{l -
l_1}{l - l_2}> 0$. On rearranging~\eqref{e:keven}, we thus obtain
the inequality~\eqref{e:kodd} for $k$ even as well.
Using~\eqref{e:kodd}, we may then choose suitable constants $\alpha$
and $\beta$, and conclude that $g \preceq G$.
\end{proof}

By the proof of the \fullref{thm:functor}, it follows that for any $v > 0$
there are uncountably many commensurability classes of nonuniform
lattices in $\Aut(X)$ of covolume $v$. This completes the proof
of~\eqref{i:nonunifcount}.

%%%%%%%%%%%%%%%%%%%%%%%%%%%%%%%%%%%%%%%%%%%%%%%%%%%%%%%%%%%%%%%%%%%%%%%%%%%%%
\subsection{Towers}
%%%%%%%%%%%%%%%%%%%%%%%%%%%%%%%%%%%%%%%%%%%%%%%%%%%%%%%%%%%%%%%%%%%%%%%%%%%%%

We construct towers in $\Aut(X)$ using full subcomplexes of
subgroups, as explained in \fullref{ss:towers}.  Let $T$ be a
biregular tree, with at least one valence greater than $2$. By the
proof of the \fullref{thm:functor}, it suffices to find a tower of tree
lattices in $\Aut(T)$, constructed as a sequence of full subgraphs
of subgroups, over a graph $A$ whose vertices may be $2$--coloured.

%%%%%%%%%%%%%%%%%%%%%%%%%%%%%%%%%%%%%%%%%%%%%%%%%%%%%%%%%%%%%%%%%%%%%%%%%%%%
\subsubsection[Proof of~\ref{i:uniftower}]{Proof
of~\eqref{i:uniftower}}
%%%%%%%%%%%%%%%%%%%%%%%%%%%%%%%%%%%%%%%%%%%%%%%%%%%%%%%%%%%%%%%%%%%%%%%%%%%%

Let $m$ and $q$ be integers $\geq 2$.  Bass and
Kulkarni~\cite[Proposition~7.15]{bk1:utl}, give an example of a tower
of uniform tree lattices acting on the $(m+1,q)$--biregular tree.
The vertices of the graph $A$ in this example may be $2$--coloured.

%%%%%%%%%%%%%%%%%%%%%%%%%%%%%%%%%%%%%%%%%%%%%%%%%%%%%%%%%%%%%%%%%%%%%%%%%%%%
\subsubsection[Proof of~\ref{i:nonuniftower}]{Proof
of~\eqref{i:nonuniftower}}
%%%%%%%%%%%%%%%%%%%%%%%%%%%%%%%%%%%%%%%%%%%%%%%%%%%%%%%%%%%%%%%%%%%%%%%%%%%%

An \emph{indexing} on a graph $A$ is a map $I$ from the edges of
$A$ to the positive integers.  A graph of groups $\mathbb{A}=(A,\cA)$ over $A$
induces the indexing $I(e)=[\cA_{i(e)}:\alpha_e(\cA_e)]$. We say that an
indexed graph $(A,I)$ admits a tower of lattices if there exists a sequence of
full subgraphs of subgroups over $A$ which induces a tower of tree lattices,
such that each graph of groups in the sequence induces the indexing $I$.

Carbone and Rosenberg~\cite[Theorem~5.3]{cr1:it} give a sufficient
condition for $(A,I)$ to admit a tower of nonuniform tree lattices.
It is easy to construct examples of $(A,I)$ which satisfy this
condition, so that the tower of tree lattices acts on the regular or
biregular tree, and the vertices of $A$ are $2$--colourable.

%%%%%%%%%%%%%%%%%%%%%%%%%%%%%%%%%%%%%%%%%%%%%%%%%%%%%%%%%%%%%%%%%%%%%%%%%%%%
\subsubsection[Proof of~\ref{i:gs}]{Proof
of~\eqref{i:gs}}
%%%%%%%%%%%%%%%%%%%%%%%%%%%%%%%%%%%%%%%%%%%%%%%%%%%%%%%%%%%%%%%%%%%%%%%%%%%%

We apply \cite[Example~7.13]{bk1:utl}. The graph $A$ in this
construction is a single edge, so its vertices are $2$--colourable.
The universal covering tree is the $(mp,q)$--biregular tree, where
$m$, $p$ and $q$ are positive integers $\geq 2$.  Putting $q_i = mp$
which is composite, there exists a tower of uniform lattices in
$\Aut(X)$ with quotient a single chamber.

\bibliographystyle{gtart}
\bibliography{link}

\begin{thebibliography}{}
\providecommand\bibmarginpar{\leavevmode\marginpar}
\def\urlstyle#1{{\tt #1}}

\bibitem{b1:ctgg}
\textbf{H Bass}, \href{http://dx.doi.org/10.1016/0022-4049(93)90085-8}
  {\emph{Covering theory for graphs of groups}}, J. Pure Appl. Algebra 89
  (1993) 3--47 \xox{MR}{1239551}

\bibitem{bk1:utl}
\textbf{H Bass}, \textbf{R Kulkarni},
  \href{http://links.jstor.org/sici?sici=0894-0347(199010)3:4%3C843:UTL%3E2.0.%
CO%3B2-B} {\emph{Uniform tree lattices}}, J. Amer. Math. Soc. 3 (1990) 843--902
  \xox{MR}{1065928}

\bibitem{bl1:tl}
\textbf{H Bass}, \textbf{A Lubotzky}, \emph{Tree lattices}, Progress in
  Mathematics 176, Birkh\"auser, Boston (2001) \xox{MR}{1794898}

\bibitem{b1:ds}
\textbf{A Borel}, \emph{On the set of discrete subgroups of bounded covolume in
  a semisimple group}, Proc. Indian Acad. Sci. Math. Sci. 97 (1987) 45--52
  (1988) \xox{MR}{983603}

\bibitem{b1:ih}
\textbf{M Bourdon}, \href{http://dx.doi.org/10.1007/PL00001619}
  {\emph{Immeubles hyperboliques, dimension conforme et rigidit\'e de
  {M}ostow}}, Geom. Funct. Anal. 7 (1997) 245--268 \xox{MR}{1445387}

\bibitem{bh1:ms}
\textbf{M\,R Bridson}, \textbf{A Haefliger}, \emph{Metric spaces of
  non-positive curvature}, Grundlehren series 319, Springer, Berlin (1999)
  \xox{MR}{1744486}

\bibitem{bm1:lpt}
\textbf{M Burger}, \textbf{S Mozes},
  \href{http://www.numdam.org/item?id=PMIHES_2000__92__151_0} {\emph{Lattices
  in product of trees}}, Inst. Hautes \'Etudes Sci. Publ. Math.  (2000)
  151--194 (2001) \xox{MR}{1839489}

\bibitem{cr1:it}
\textbf{L Carbone}, \textbf{G Rosenberg}, \emph{Infinite towers of tree
  lattices}, Math. Res. Lett. 8 (2001) 469--477 \xox{MR}{1849263}

\bibitem{fh1:ci}
\textbf{B Farb}, \textbf{G\,C Hruska}, \emph{Commensurability invariants for
  nonuniform tree lattices}, Israel J. Math. 152 (2006) 125--142
  \xox{MR}{2214456}

\bibitem{gp1:ih}
\textbf{D Gaboriau}, \textbf{F Paulin},
  \href{http://dx.doi.org/10.1023/A:1013168623727} {\emph{Sur les immeubles
  hyperboliques}}, Geom. Dedicata 88 (2001) 153--197 \xox{MR}{1877215}

\bibitem{g1:gs}
\textbf{Y Glasner}, \href{http://dx.doi.org/10.1016/S0021-8693(03)00530-1}
  {\emph{A two-dimensional version of the Goldschmidt--Sims conjecture}}, J.
  Algebra 269 (2003) 381--401 \xox{MR}{2015283}

\bibitem{g1:atg}
\textbf{D\,M Goldschmidt},
  \href{http://links.jstor.org/sici?sici=0003-486X(198003)2:111:2%3C377:AOTG%3%
E2.0.CO%3B2-U} {\emph{Automorphisms of trivalent graphs}}, Ann. of Math. $(2)$
  111 (1980) 377--406 \xox{MR}{569075}

\bibitem{hp1:cai}
\textbf{F Haglund}, \textbf{F Paulin},
  \href{http://dx.doi.org/10.1007/s00208-002-0373-x} {\emph{Constructions
  arborescentes d'immeubles}}, Math. Ann. 325 (2003) 137--164 \xox{MR}{1957268}

\bibitem{js1:hcg}
\textbf{T Januszkiewicz}, \textbf{J {\'S}wi{\c{a}}tkowski},
  \href{http://dx.doi.org/10.1007/s00014-003-0763-z} {\emph{Hyperbolic Coxeter
  groups of large dimension}}, Comment. Math. Helv. 78 (2003) 555--583
  \xox{MR}{1998394}

\bibitem{l1:tl}
\textbf{A Lubotzky}, \emph{Tree-lattices and lattices in Lie groups}, from:
  ``Combinatorial and geometric group theory (Edinburgh, 1993)'', London Math.
  Soc. Lecture Note Ser. 204, Cambridge Univ. Press (1995)  217--232
  \xox{MR}{1320284}

\bibitem{pv1:rrg}
\textbf{L Potyagailo}, \textbf{E Vinberg}, \emph{On right-angled reflection
  groups in hyperbolic spaces}, Comment. Math. Helv. 80 (2005) 63--73
  \xox{MR}{2130566}

\bibitem{ros1:tctl}
\textbf{G\,E Rosenberg}, \emph{Towers and Covolumes of Tree Lattices}, PhD
  thesis, Columbia University (2001)

\bibitem{s:cgd}
\textbf{J-P Serre}, \emph{Cohomologie des groupes discrets}, from: ``Prospects
  in mathematics (Proc. Sympos. Princeton 1970)'', Ann. of Math. Studies 70,
  Princeton Univ. Press (1971)  77--169 \xox{MR}{0385006}

\bibitem{s1:t}
\textbf{J-P Serre}, \emph{Trees}, Springer, Berlin (1980) \xox{MR}{607504}

\bibitem{at1:cul}
\textbf{A Thomas}, \emph{Covolumes of uniform lattices acting on polyhedral
  complexes}, Bull. London Math. Soc. to appear

\bibitem{vs1:g}
\textbf{{\`E}\,B Vinberg}, \textbf{O\,V Shvartsman}, \emph{Discrete groups of
  motions of spaces of constant curvature}, from: ``Geometry II: Spaces of
  Constant Curvature'', Encyclopaedia Math. Sci. 29, Springer, Berlin (1993)
  139--248 \xox{MR}{1254933}

\end{thebibliography}

\end{document}